\documentclass[12pt]{article}
\usepackage{graphicx}
\usepackage{amsfonts}
\usepackage{latexsym}
\usepackage{amsmath}

\makeatletter
\@addtoreset{figure}{section}
\def\thefigure{\thesection.\@arabic\c@figure}
\@addtoreset{table}{section}
\def\thetable{\thesection.\@arabic\c@table}

\def\@sect#1#2#3#4#5#6[#7]#8{\ifnum #2>\c@secnumdepth
     \def\@svsec{}\else
     \refstepcounter{#1}\edef\@svsec{\csname the#1\endcsname.\hskip .75em
}\fi
     \@tempskipa #5\relax
      \ifdim \@tempskipa>\z@
        \begingroup #6\relax
          \@hangfrom{\hskip #3\relax\@svsec}{\interlinepenalty \@M #8\par}%
        \endgroup
       \csname #1mark\endcsname{#7}\addcontentsline
         {toc}{#1}{\ifnum #2>\c@secnumdepth \else
                      \protect\numberline{\csname the#1\endcsname}\fi
                    #7}\else
        \def\@svsechd{#6\hskip #3\@svsec #8\csname #1mark\endcsname
                      {#7}\addcontentsline
                           {toc}{#1}{\ifnum #2>\c@secnumdepth \else
                             \protect\numberline{\csname the#1\endcsname}\fi
                       #7}}\fi
     \@xsect{#5}}
\def\@begintheorem#1#2{\it \trivlist \item[\hskip \labelsep{\bf #1\ #2.}]}
\def\section{\@startsection {section}{1}{\z@}{-3.5ex plus -1ex minus
 -.2ex}{2.3ex plus .2ex}{\normalsize\bf}}

\begin{document}

\title{Dense Packings of Congruent Circles \\
in Rectangles with a Variable Aspect Ratio}
\date{} 
\maketitle

\begin{center}
\author{Boris D. Lubachevsky  \ \ \ \ \ \ \ \ \ \ \ \ \ \ \ \ \ \ \ \ \ \ \ \ \ Ronald  Graham\\
{\em bdl@bell-labs.com \ \ \ \ \ \ \ \ \ \ \ \ \ \ \ \ \ \ \ \ \ \ \ \ \ \ \ graham@ucsd.edu}\\
Bell Laboratories \ \ \ \ \ \ \ \ \ \ \ \ \ \ \ \ \ \ \ \ \ \ University of California \\
600 Mountain Avenue  \ \ \ \ \ \ \ \ \ \ \ \ \ \ \ \ \ \ \ \ \ \ \ \ \ \ \ \ \ at San Diego \\
Murray Hill, New Jersey \ \ \ \ \ \ \ \ \ \ \ \ \ \ \ \ \ \ La Jolla, California }
\end{center}

\setlength{\baselineskip}{0.995\baselineskip}
\normalsize
\vspace{0.5\baselineskip}
\vspace{1.5\baselineskip}

\begin{abstract}
We use computational experiments 
to find the rectangles of minimum area into which a given number
$n$ of non-overlapping congruent circles can be packed. No
assumption is made on the shape of the rectangles. 
Most of the packings
found have the usual regular square or hexagonal
pattern. However,
for 1495 values of $n$ in the tested
range $n \le 5000$, specifically, for
$n = 49, 61, 79, 97, 107,... 4999$,
we prove that the optimum cannot possibly be achieved
by such regular arrangements.
The evidence suggests that the limiting height-to-width ratio
of rectangles containing an optimal hexagonal packing of circles
tends to $2 - \sqrt {3}$ as $n \rightarrow \infty$,
if the limit exists.



\end{abstract}
\section{Introduction}
\hspace*{\parindent} 
Consider the task of finding
the smallest area rectangular region
that encloses a given number $n$ of
circular disks of equal diameter.
The circles must not
overlap with each other or extend outside the rectangle.
The 
aspect ratio of the rectangle,
i.e., the ratio of its height to width,
is variable and 
subject to the area-minimizing choice
as well as the positions of the circles
inside the rectangle.

\begin{figure}
\centering
\includegraphics*[width=6.5in,height=3.9in]{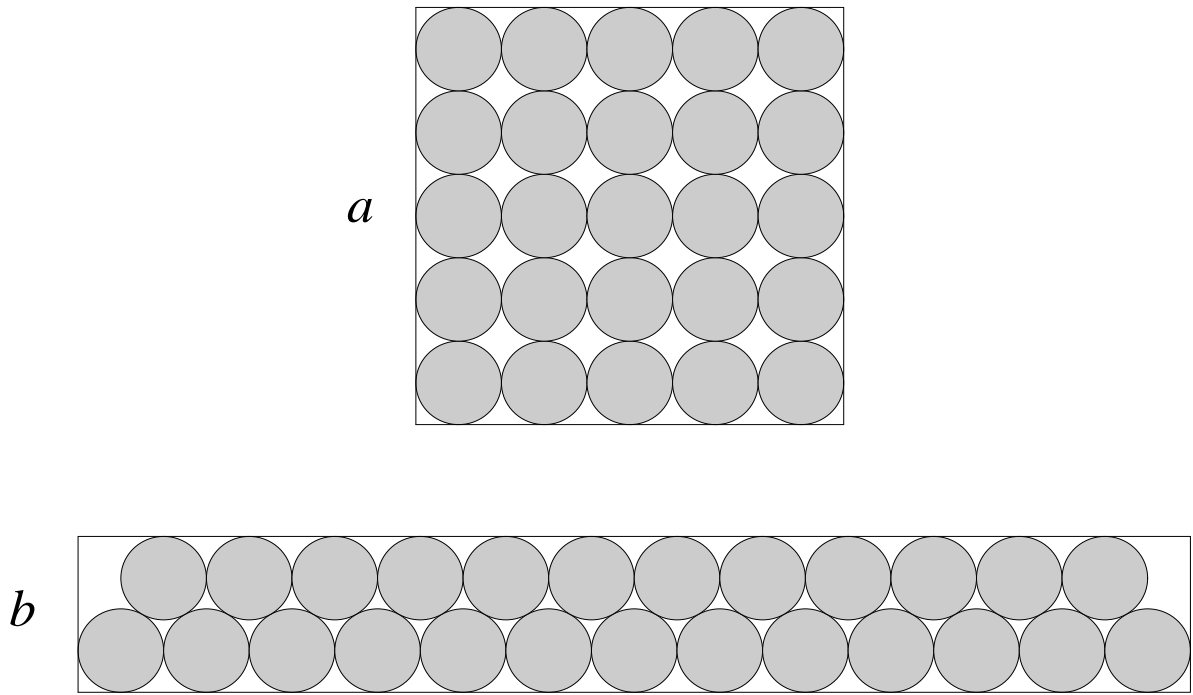}
\caption{The best packings found for 25 circles $a$) in a square,
and, $b$) in a rectangle with variable aspect ratio. 
}
\label{fig:25}
\end{figure}

Packing circles 
in a square has been the subject of
many investigations \cite{GL}, \cite{NO1}, \cite{NO2}, \cite{NO3}.
Because the aspect ratio is not fixed in our present problem,
the solutions are typically different from the dense
packings in a square.
For example, the density $\pi/4 = 0.785... $ of the proved 
optimum 
(see \cite{NO3})
for packing 25 congruent circles in a square 
(see Figure~\ref{fig:25}a)
can be increased to $25\pi/(26(2 + \sqrt 3 )) = 0.809...$,
if we 
let the rectangle assume its best aspect ratio
(see Figure~\ref{fig:25}b).
Our experiments offer
only three values of $n$ for which 
the dense packing in a square is also a solution to our
present problem: $n = 4, 9$, and of course, $n=1$.

\begin{figure}
\centering
\includegraphics*[width=6.5in,height=1.5in]{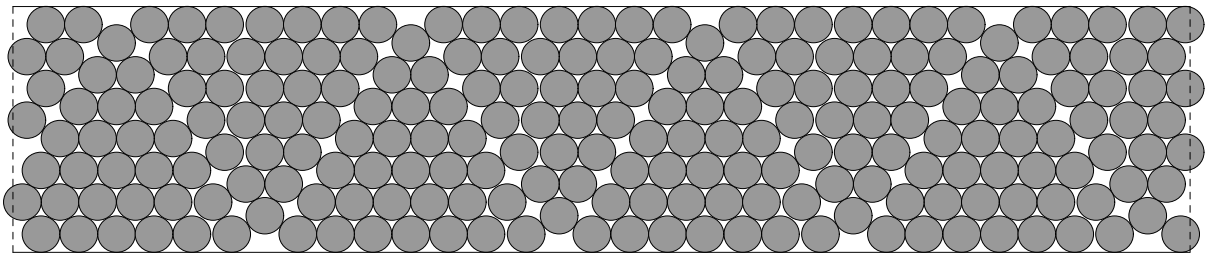}
\caption{A 224-circle fragment of the packing with the highest density found
in a long rectangle with a fixed aspect ratio.
}
\label{fig:224}
\end{figure}

Similarly, 
dense packings in long rectangles with a fixed aspect ratio
usually do not yield solutions for our problem.
According to a long-standing conjecture
attributed to Molnar (see \cite{Furedi}), 
dense packings in long fixed rectangles
tend to form 
periodic up and down 
alternating 
triangular ``teeth,''
as in the example shown in Figure~\ref{fig:224}.
We can usually increase the density of such packings
by slightly changing the aspect ratio
of the rectangle.

The present problem has been occasionally mentioned
among packing problems ;
for example, Web enthusiasts have recently began discussing it.
However, 
the problem was considered as early as in 1970 in \cite{Ruda},
where the optimum packings were determined for all $n \le 8$
and conjectured for $9 \le n \le 12$.
We learned about Ruda's results 
after we performed our computations and made our conjectures.
Fortunately,
our conjectures restricted to the values $n \le 12$ agree with
the results and conjectures in \cite{Ruda}!

%
It appears that the complexity of the proofs of
optimality in the packing-circles-in-a-square problem
increases exponentially with $n$;
computer generated/assisted proofs of optimality
of such packings do not
reach $n = 30$ (see \cite{NO3} ;
proofs that do not utilize computers
have not been found except for quite small values
of $n$),
while conjectures extend to $n = 50$ and even much larger values
(see \cite{NOR}).

Similarly, it seems to be difficult for the present problem
to prove optimality for the configurations we find.
This paper describes an experimental approach,
where good packings are obtained with the help of a computer.
Some of these packings
hopefully will be proved optimal
in the future. At present, most of the statements 
about the packings made
in this paper are only conjectures,
except for a few which we explicitly claim as proven.
For instance, we prove that the best packings
found are better than any other in their class,
e.g., that the packing of 79 circles with a monovacancy
(see Figure~\ref{fig:79}a),
while non-optimal,
is better than any hexagonal or square grid packing
of 79 circles without a monovacancy.
Also, we prove that $n = 11$ is the smallest $n$ for which a hexagonal
packing as in Figure~\ref{fig:11} is better than any square grid packing
of $n$ circles.
\section{A priori expectations and questions about best packings}\label{sec:hexandsq}
\begin{figure}
\centering
\includegraphics*[width=1.6in,height=1.13in]{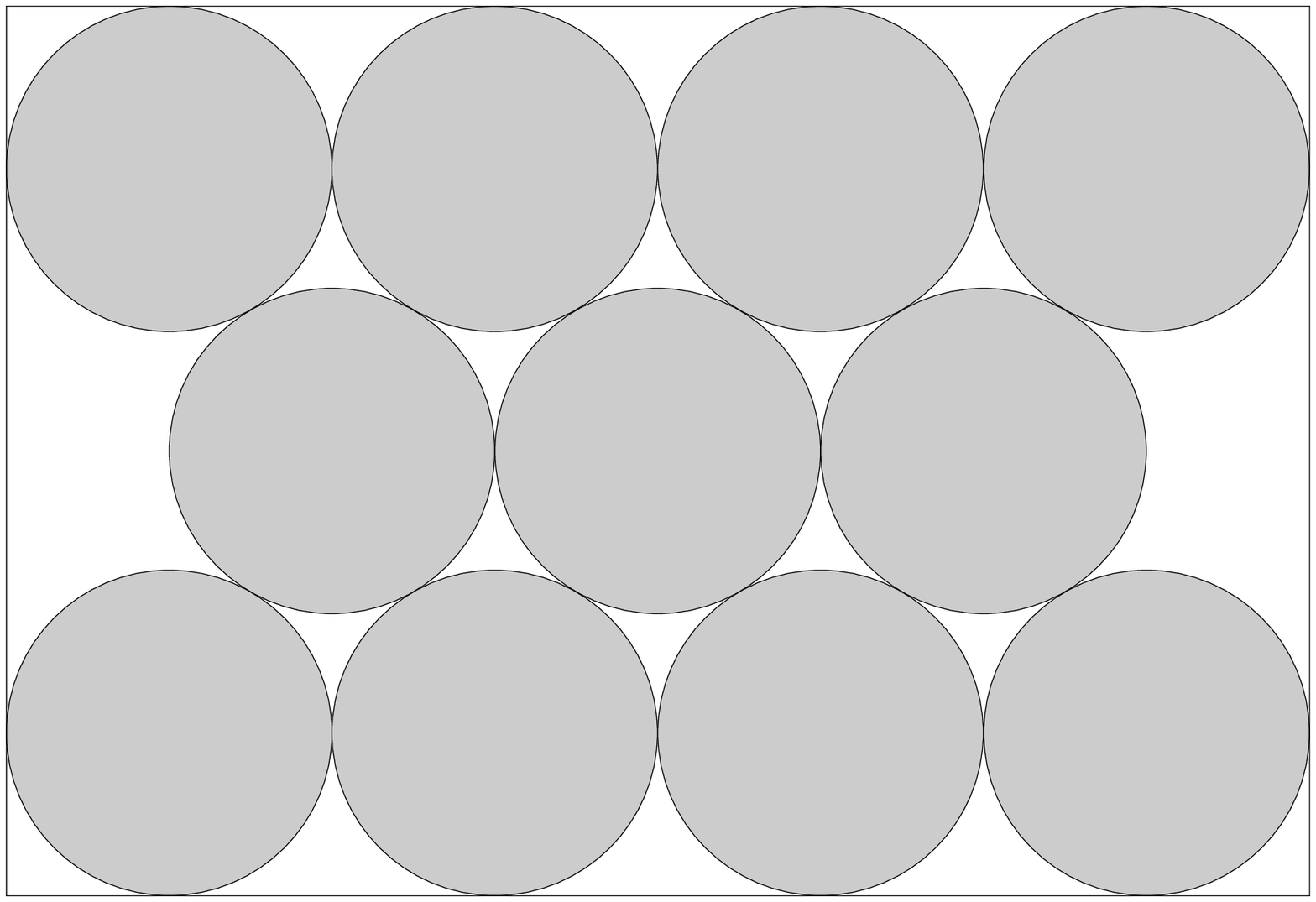}
\caption{The best packing found for 11 circles in a rectangle.
}
\label{fig:11}
\end{figure}

\begin{figure}
\centering
\includegraphics*[width=7.1in,height=4.5in]{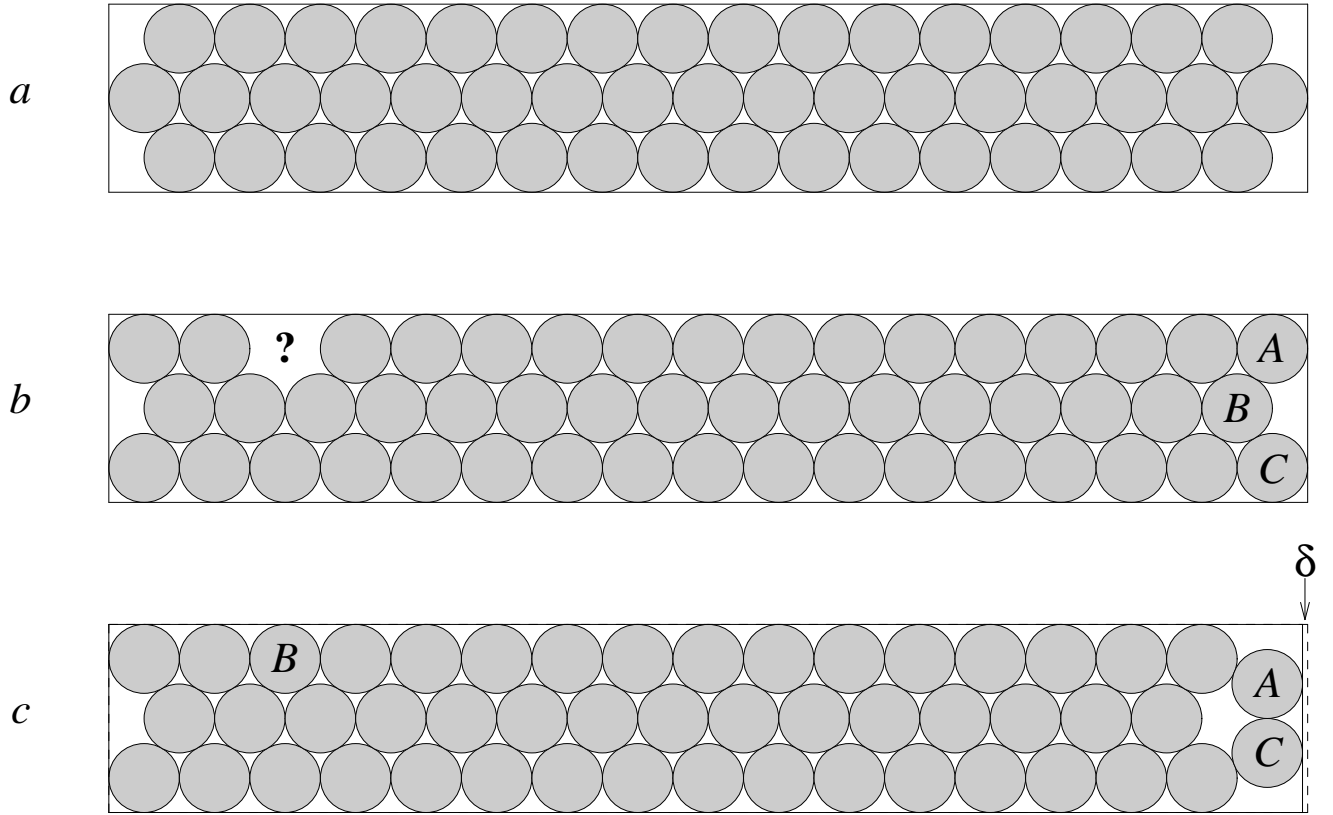}
\caption{Packings of 49 circles in a rectangle: 
$a$) the best in the class of hexagonal packings,
$b$) a best in the class of hexagonal packings with monovacancies
(one of 17 equally dense
packings
with
the hole),
$c$) the best we found. 
}
\label{fig:49}
\end{figure}

\begin{figure}
\centering
\includegraphics*[width=6.5in,height=3.9in]{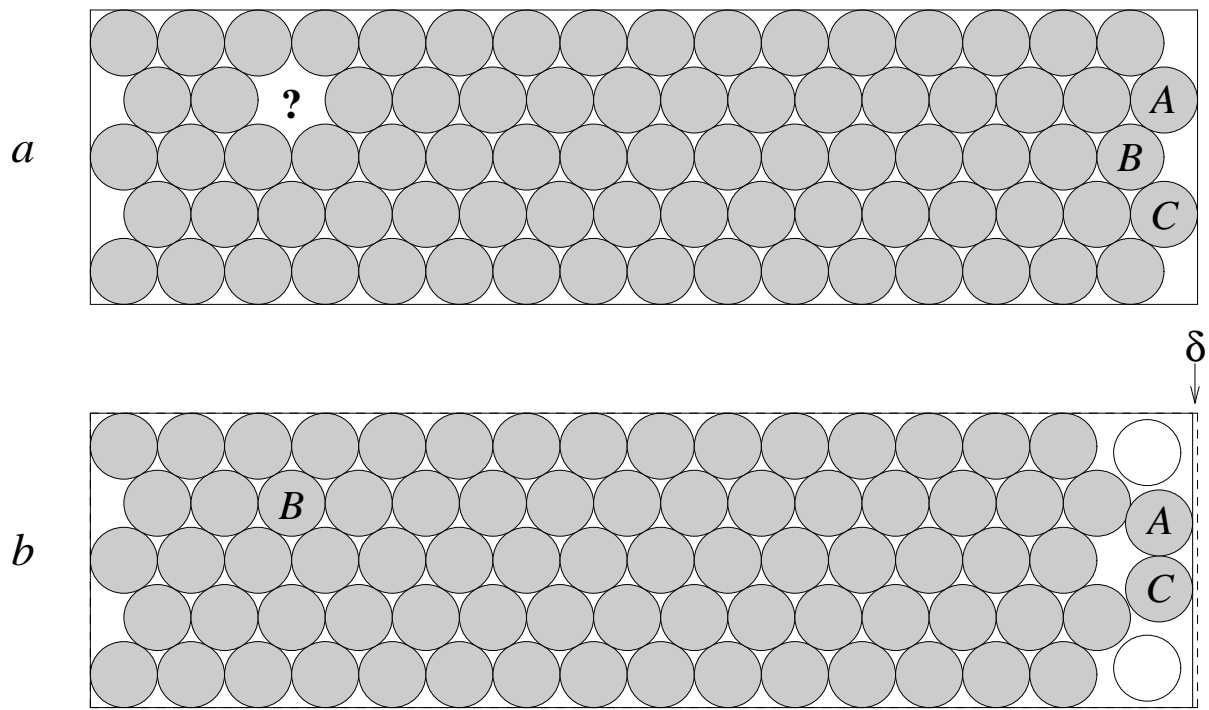}
\caption{Packings of 79 circles in a rectangle:
$a$) a best in the class of hexagonal packings with possible monovacancies,
and, $b$) the best we found.
}
\label{fig:79}
\end{figure}

\hspace*{\parindent} 
Since it is well known that the hexagonal arrangement of 
congruent circles
in the infinite plane has the highest possible density,
(see \cite{Th}, \cite{Ft}, 
\cite{FG},\cite{Oler}), 
before plunging into our experiments,
we expected
to obtain good finite packings by
``carving'' rectangular subsets out of this infinite packing.
It will be useful to classify here such finite arrangements.
For that we will discuss packings in
Figures  \ref{fig:25}b, \ref{fig:11}, 
\ref{fig:49}a,
\ref{fig:49}b, 
and \ref{fig:79}a, irrespective of their
optimality (which will be discussed in the following sections).
For the purpose of exposition,
we will pretend in this section that there are no holes in
the arrangements in Figures~\ref{fig:49}b and \ref{fig:79}a,
that is, the question mark in the figures is covered with 
an additional circle of the common radius.

\begin{figure}
\centering
\includegraphics*[width=6.5in,height=2.7in]{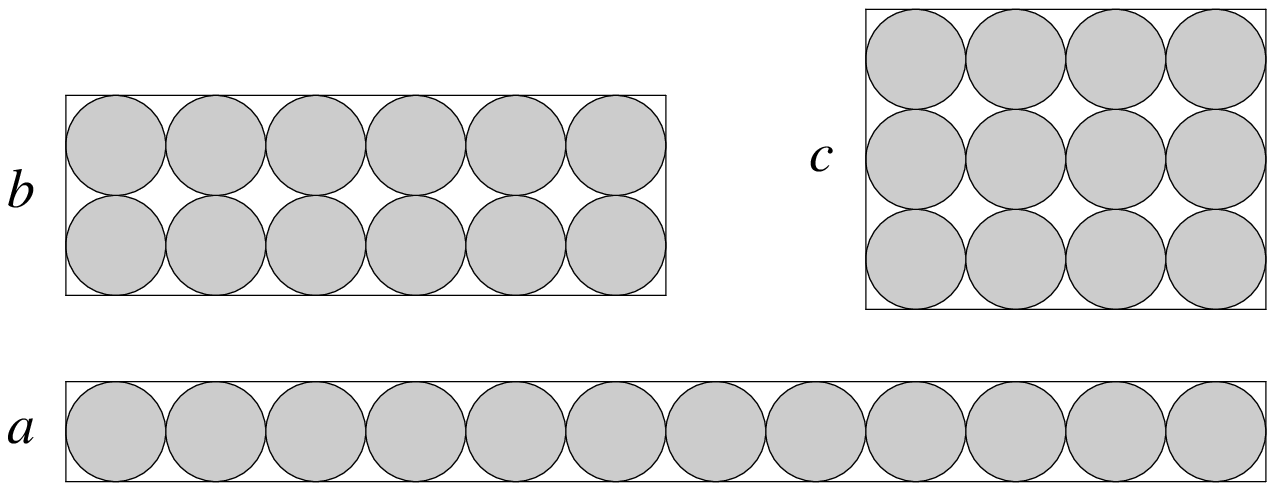}
\caption{Best packings found for 12 circles in a rectangle.
}
\label{fig:12}
\end{figure}

Under such an assumption,
the arrangement 
in Figure~\ref{fig:79}a
consists of $h=5$ alternating rows
with $w=16$ circles in each,
a total of $n = w \times h =80$ circles.
Arrangements as in Figures~\ref{fig:25}b, ~\ref{fig:11}, \ref{fig:49}a,
and 
\ref{fig:49}b,
are obtained by depleting such
$w \times h$ arrays of some circles along one side of the array.
The ways of performing the depletion depend on the parity of $h$.

If $h$ is even, then the depletion can be done in
only one way by removing $h/2$ circles.
Figure~\ref{fig:25}b presents an example with $h = 2$.

The case of an odd $h$ presents two possibilities:
we can remove $\lfloor h/2 \rfloor$ circles on a side
as in the example of Figure~\ref{fig:49}b,
or $\lfloor h/2 \rfloor + 1$ circles on a side
as in the example of Figure~\ref{fig:49}a.

Alternatively, one can consider rectangular square grid
arrangements as candidates for the best packings,
e.g., see those in Figure~\ref{fig:12}.
The density of a square grid arrangement
is fixed at $\pi /4$ independent of $n$. 
When both sides $h$ and $w$ of the carved rectangle tend to infinity,
the hexagonal packing density tends to $\pi /(2 \sqrt 3 )$ which is
the density of the hexagonal packing in the infinite plane.
Since $\pi / (2 \sqrt 3 ) > \pi /4$,
one naturally wonders
for which $n$ the best hexagonal grid arrangement becomes
better than the square grid arrangement.

A natural question is whether or not
both arrangements exhaust all possibilities for the best packing.
In other words, does there exist
an optimal packing of $n$ disks in a rectangle
such that it {\em cannot} be represented as being carved out 
by a rectangle from
either square grid or hexagonal grid
packing of the infinite plane?

\section{Results of compactor simulations}\label{sec:simulation}
\hspace*{\parindent}
To tackle the problem by computer, we developed
a ``compactor'' simulation algorithm.
The simulation begins by starting with a random initial
configuration with $n$ circles lying inside a (large) rectangle
without circle-circle overlaps. The starting configuration
is feasible but is usually rather sparse.
Then the computer imitates a ``compactor'' with
each side of the rectangle pressing against the circles,
so that the circles are being forced
towards each other until they ``jam.''
Possible circle-circle or circle-boundary conflicts
are resolved using a simulation of a hard collision so
that no overlaps occur during the process.

The simulation for a particular $n$ is repeated many times,
with different starting circle configurations.
If the final density in a run is larger than
the record achieved thus far, it replaces this record.
Eventually in this process, the record stops improving up to
the level of accuracy induced by the double precision accuracy
of the computer. The resulting packing now becomes a candidate
for the optimal packing for this value of $n$.

In this way we found that the square grid pattern 
with density $\pi /4$
supplies
the optimum for $n=1,...,10$, and that
$n=11$ is the smallest number of circles
for which a density better than $\pi /4$
can be reached; 
the density of this packing is $11 \pi /(16(1+ \sqrt 3 )= 0.790558..$
The corresponding conjectured
optimum pattern for $n = 11$ is shown in
Figure~\ref{fig:11}.
The pattern is hexagonal.
For $n=12$ and 13 the square grid pattern
briefly takes over as the experimental optimum again
(see Figure~\ref{fig:12} for the case of $n = 12$).
Remarkably, Ruda \cite{Ruda} {\em proves} that for $n \le 8$, the 
square-grid packings are optimal. 
His conjectures for $9 \le n \le 12$ also agree
with our findings.

Our simulation results show that
for $n \ge 14$, the best densities are larger than $\pi/4$.
This statement is also easy to prove, considering
examples of (not necessarily optimal)
two-row hexagonal packings 
like that in Figure~\ref{fig:25}b for odd $n = 2m + 1 > 14$,
or the ones with equal length rows for even $n = 2m \ge 14$.
Assuming the circles have radius 1,
the following inequalities
\begin{equation}
\label{odd}
  2(m+1)(2 + \sqrt 3 ) < 4(2m + 1)
\end{equation}
for $n = 2m + 1$
and
\begin{equation}
\label{even}
  (2m+1)(2 + \sqrt 3 ) < 8m
\end{equation}
for $n = 2m$ have to be satisfied.
The left-hand sides in \eqref{odd} and \eqref{even}
are the areas of the enclosing rectangles for the
two-row hexagonal packings, and the right-hand sides
are the areas for the corresponding square grid packings.
It is easily seen that all integers $m \ge 7$ satisfy either inequality.
These correspond to all integer values of $n \ge 14$.

\begin{figure}
\centering
\includegraphics*[width=4.0in,height=2.89in]{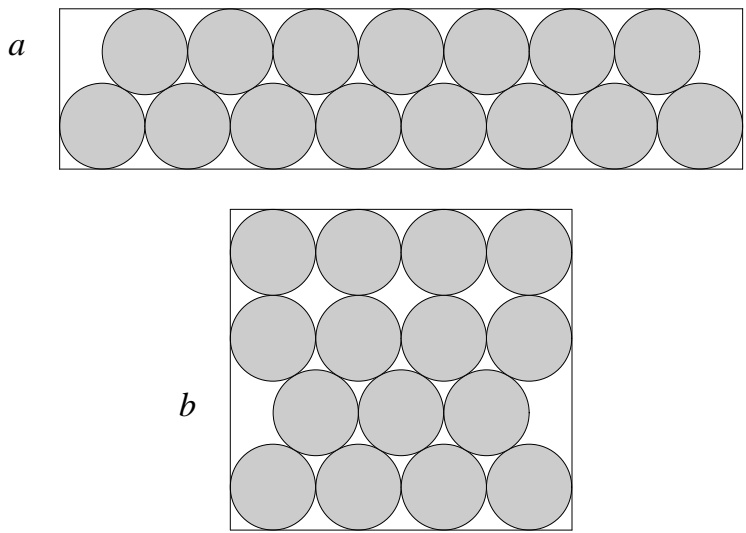}
\caption{Two equally dense
best packings found for 15 circles in a rectangle.
}
\label{fig:15}
\end{figure}

A surprise awaited us at the value $n = 15$.
We found two different equally dense packings.
One, shown in Figure~\ref{fig:15}a, is of the expected hexagonal type.
However, the other, shown in Figure~\ref{fig:15}b, is not, nor can it be
carved out of the square grid.
It is easy to verify that
two equally dense packings $a$ and $b$ of these types also exist
for any $n$ of the form $n = 15 + 4k$, $k=1,2,3, ...$.
Packing $a$, such as
that in Figure~\ref{fig:15}a,
has two alternating rows, with one row one circle shorter
than the other, and with the longer row consisting 
of $w = 8 + 2k$ circles.
Packing $b$, such as
that in Figure~\ref{fig:15}b,
has 4 rows,
the longest row having $w = 4 + k$ circles,
three bottom rows alternate, the middle one having $w - 1$ circles,
and the 4th row of $w$ circles stacked straight on the top
(or on the bottom; these would produce the same configuration).
Our experiments suggest that
for $k = 1$ and 4, i.e., for $n = 19$ and 31,
these pairs of configurations might also be optimal;
however, they are provably not optimal for other values of $k > 0$.
The packings of pattern $b$ for $n=15$, 19, and 31, if
they are proved to be optimal, answer positively
the second question posed in Section~\ref{sec:hexandsq}.

Thus, the optimal container for 15 bottles 
apparently can have two different shapes,
as can the optimal containers for 19 and 31 bottles!
In an obvious way, 
more than one optimum container shape also exists
for square grid packings of $n$ circles if $n$ is not a prime.
See Figure~\ref{fig:12} for the example of $n = 12$;
there are three different rectangles which are equally good
and probably optimal for $n = 12$.
Apparently, three is the maximum number of rectangular shapes,
any number $n$ of circles can optimally fit in.
It seems that
for 
$n = 4, 6, 8, 9, 10, 15, 19$ and 31,
there are exactly two shapes.
For all other $n$ tested,
only one best aspect ratio was found experimentally.

Another packing surprise awaited us at $n = 49$.
Figures~\ref{fig:49}a and \ref{fig:49}b show two among several
equivalent
packings achieved in our simulation experiments.
The configuration $a$ is hexagonal, while
configuration $b$ contains a monovacancy,
that is, a hole than can accommodate exactly one circle.
Monovacancies in hexagonal arrays of congruent circles appear
often in the simulation experiments
for $14 \le n < 49$ but only in packings of
inferior quality, 
so that a better 
quality hexagonal packing
without monovacancies
could always be found.
No higher density hexagonal
packing without monovacancies
was found for $n = 49$.

Such large $n$ already present substantial
difficulties for our simulation procedure.
The procedure
failed to produce a packing which is better
than those in Figures~\ref{fig:49}a and \ref{fig:49}b.
However, it
is easy to prove that the density of a hexagonal packing in a rectangle
with a monovacancy can be increased.
For example, to
improve the packing in Figure~\ref{fig:49}b
we 
relocate circle $B$ 
into the vacancy 
and then rearrange circles $A$ and $C$
along the side.
This reduces the width of the rectangle
by a small but positive $\delta$ as shown
in Figure~\ref{fig:49}c,
where $\delta = 2 - \sqrt { 2 \sqrt 3 } = 0.13879..$
of the circle radius.
The density of $c$ is $49 \pi/(2(1+ \sqrt 3 )(34- \delta)) = 0.83200266..$.

It is not known whether or not
the resulting configuration $c$ can be further
improved. However, using the exhaustive search
we show (see next section)
that the optimum packing for
49 circles in a rectangle,
whatever it might be,
cannot be purely
hexagonal. 
The value $n=49$ is apparently the smallest
one for which neither 
square-grid nor perfect hexagonal pattern
delivers the optimum.

\section{Exhaustive search}\label{sec:exhaust}
\hspace*{\parindent}
The simulation method becomes progressively slower
for increasing $n$.
However, by exercising the simulation 
for smaller $n$ we are able to refine the 
idea of the class of packings
which might deliver the optimum.
We have chosen the class that consists of 
square grid and hexagonal packings and their hybrids;
we also allow for 
monovacancies in configurations of the class,
although these configurations are never optimal.
We have extended our computing experiments to larger values of $n$
using exhaustive search.
For a given $n$,
the method
simply evaluates 
each candidate in the class by
computing the area of the rectangle
and selects
the one with the smallest area.
Note that
for each $n$ in this class, there are only a finite number
of candidates.

\begin{figure}
\centering
\includegraphics*[width=2.45in,height=3.09in]{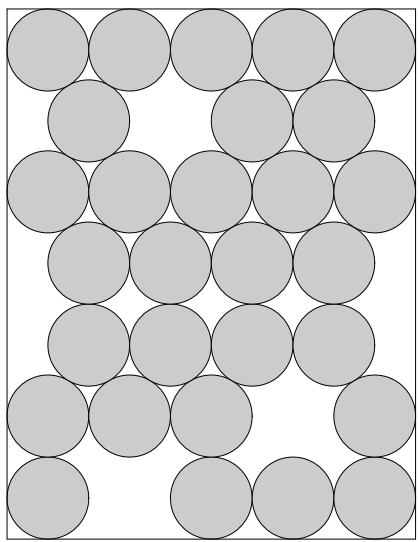}
\caption{A ``general case'' member of the class of packings
among which we search for the optimum; 
here $w = 5$, $h = 5$, $h_- = 2$, $s = 2$, $s_{-} = 1$, and $d = 3$.
}
\label{fig:class}
\end{figure}

A general packing in the class
(see Figure~\ref{fig:class}),
consists of $h + s$ rows and has $d$ monovacancies.
The $h$ rows are hexagonally alternating, and
the $s$ rows are stacked 
directly on top of the previous
row as in square grid packings.
The longest row consists of $w$ circles.
Assuming all monovacancies are filled with circles,
among the $h$ hexagonal rows, 
$h_{-}$ rows consist of $w-1$ circles each,
the remaining $h - h_{-}$ rows consist of $w$ circles each,
and
among the $s$ square grid rows, 
$s_{-}$ rows consist of $w-1$ circles each,
and the remaining $s - s_{-}$ rows consist of $w$ circles each.
The number of circles in the configuration is
\begin{equation}
\label{ncircles}
n = w(h + s) - h_{-} - s_{-} - d
\end{equation}

Table~\ref{tab:1t53} lists the best packings found
of $n$ circles in a rectangle with variable aspect ratio
for $n$ in the range $1 \le n \le 53$,
and Table~\ref{tab:54t213} continues this list
for $54 \le n \le 213$.
The packings in Table~\ref{tab:1t53} are obtained
by the simulation procedure described in Section~\ref{sec:simulation}
and verified by the exhaustive search.
Most packings in Table~\ref{tab:54t213} were generated
only by the exhaustive search.
An entry in either table consists of the $n$, followed by
the set of integers which represent
the parameters of the packing structure as
explained above:
\\
$~~~~w$, the number of circles in the longest row,
\\
$~~~~h$, the number of rows arranged in a hexagonal alternating pattern,
\\
$~~~~h_{-}$, the number of rows that consist of $w-1$
circles each,
\\
$~~~~s$, the number of rows, in addition to $h$ rows, 
that are  stacked in the square grid pattern.
\\
Column $s$ is absent in Table~\ref{tab:54t213},
because, as explained in Section~\ref{sec:simulation},
no optimum packing found for $n > 31$ has  $s$ positive.
Also, 
for any $n$, 
no optimum packing was found
with 
$s_{-} > 0$; 
thus, the $s_{-}$ column is omitted. 

\begin{table}
\begin{center}
\fbox{
\begin{tabular}{r|r|r|r|r||r|r|r|r|r||r|r|r|r|r} 
$n$&$w$&$h$&$h_{-}$&$s$&$n$&$w$&$h$&$h_{-}$&$s$&$n$&$w$&$h$&$h_{-}$&$s$ \\ \hline
1  & 1 &0  & 0     & 1 &15 & 8 & 2 & 1     & 0 &33 & 7 & 5 & 2     & 0 \\
2  & 2 & 0 & 0     & 1 &   & 4 & 3 & 1     & 1 &34 & 9 & 4 & 2     & 0 \\
3  & 3 & 0 & 0     & 1 &16 & 8 & 2 & 0     & 0 &35 &12 & 3 & 1     & 0 \\
4  & 4 & 0 & 0     & 1 &17 & 6 & 3 & 1     & 0 &36 &12 & 3 & 0     & 0 \\
   & 2 & 0 & 0     & 2 &18 & 9 & 2 & 0     & 0 &37 &19 & 2 & 1     & 0 \\
5  & 5 & 0 & 0     & 1 &19 & 10& 2 & 1     & 0 &38 &13 & 3 & 1     & 0 \\
6  & 6 & 0 & 0     & 1 &   & 5 & 3 & 1     & 1 &39 &13 & 3 & 0     & 0 \\
   & 3 & 0 & 0     & 2 &20 & 7 & 3 & 1     & 0 &40 &10 & 4 & 0     & 0 \\
7  & 7 & 0 & 0     & 1 &21 & 7 & 3 & 0     & 0 &41 &14 & 3 & 1     & 0 \\
8  & 8 & 0 & 0     & 1 &22 & 11& 2 & 0     & 0 &42 &11 & 4 & 2     & 0 \\
   & 4 & 0 & 0     & 2 &23 & 8 & 3 & 1     & 0 &43 & 9 & 5 & 2     & 0 \\
9  & 9 & 0 & 0     & 1 &24 & 8 & 3 & 0     & 0 &44 &15 & 3 & 1     & 0 \\
   & 3 & 0 & 0     & 3 &25 & 13& 2 & 1     & 0 &45 &15 & 3 & 0     & 0 \\
10 &10 & 0 & 0     & 1 &26 & 9 & 3 & 1     & 0 &46 &12 & 4 & 2     & 0 \\
   & 5 & 0 & 0     & 2 &27 & 9 & 3 & 0     & 0 &47 &16 & 3 & 1     & 0 \\
11 & 6 & 2 & 1     & 0 &28 & 6 & 5 & 2     & 0 &48 &10 & 5 & 2     & 0 \\
12 &12 & 0 & 0     &12 &29 & 10& 3 & 1     & 0 &*49&17 & 3 & 2     & 0 \\
   & 6 & 0 & 0     & 2 &30 & 10& 3 & 0     & 0 &50 &17 & 3 & 1     & 0 \\
   & 4 & 0 & 0     & 3 &31 &16 & 2 & 1     & 0 &51 &17 & 3 & 0     & 0 \\
13 &13 & 0 & 0     & 1 &   & 8 & 3 & 1     & 1 &52 &13 & 4 & 0     & 0 \\
14 & 5 & 3 & 1     & 1 &32 &11 & 3 & 1     & 0 &53 &11 & 5 & 2     & 0 \\
\end{tabular}
}
\caption{Packings of $n$ circles in a rectangle 
for $1 \le n \le 53$ with $w$ circles in a row, 
$h$ hexagonal rows, 
$h_{-}$ of which
consist of $w-1$ circles each,
and $s$ square grid rows of $w$ circles each.
Except for the case $n=49$ marked with a star,
all packings are the best we could find.
}
\label{tab:1t53}
\end{center}
\end{table}

In most packings presented in these two tables,
the number of monovacancies $d = 0$,
and so the $d$ column is omitted.
A few entries where monovacancies are possible
are marked with stars.
The non-starred entries describe the best packings found.
The entries marked with stars 
cannot be optimal,
as explained in Section~\ref{sec:simulation}.

An improved packing for the marked case of 49 circles
can be obtained as described in Section~\ref{sec:simulation}
(see Figure~\ref{fig:49}).
The next marked case $n=61$ is similar;
Table~\ref{tab:1t53} lists the variant $a$ for it
with $h = 3$ and  $h_{-} = \lfloor h/2 \rfloor +1 = 2$.
The following marked case is $n=79$.
Here too we improve the packing with the hole
by relocating side circles.
$A$, $B$, $C$, $D$, and $E$
as shown in  Figure~\ref{fig:79},
where the value of improvement $\delta$ is 
the same as in the case $n=49$.

\begin{figure}
\centering
\includegraphics*[width=7in,height=3.9in]{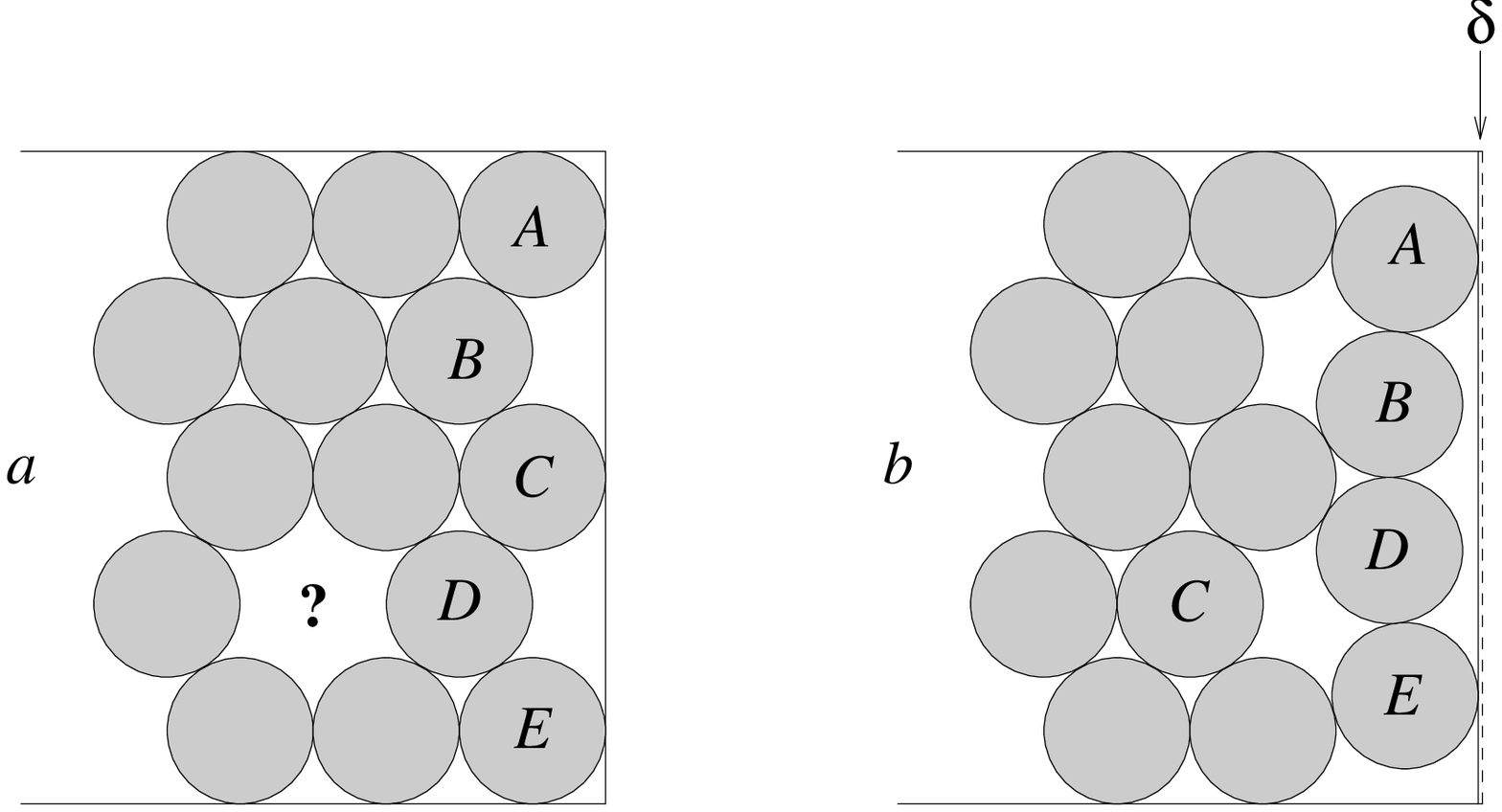}
\caption{Improving a packing with a monovacancy when
$h=5$, $h_{-} = 2$:
$a$) an original hexagonal
packing
with a monovacancy,
$b$) the transformed packing.
}
\label{fig:5h3}
\end{figure}

This method of relocation applies to
all marked cases in the tables
when $h$ is odd.
For example, for the case of $h = 5$,
when $h_{-} = 3$
the relocation can be done as illustrated
in Figure~\ref{fig:5h3}.
This applies to the cases $n = 97$, 107, and 142
in Table~\ref{tab:54t213}.
In all these cases, before
using this method, we should replace
listed in the tables
variant
that has $h_{-} = \lfloor h/2 \rfloor = 2$
with the variant
that has $h_{-} = \lfloor h/2 \rfloor + 1 = 3$.
This would produce a configuration with the pattern
of the side
as in Figure~\ref{fig:5h3}a.
We would then improve this configuration by relocating
circles $A$, $B$, $C$, $D$, and $E$
to yield a configuration as in
Figure~\ref{fig:5h3}b.
The value of the improvement here is $\delta =
2 - 0.5  \sqrt 3  -
 3^{1/4} (2 \sqrt 3 - 1)/(2 \sqrt {4 - \sqrt 3 }) = 0.05728...$
of the circle radius.
A similar method works for even $h$.

Sometimes during these improvements,
some circles become 
the so-called {\em rattlers},
i.e., they become free to move and hit their neighbors,
like the two unshaded circles in Figure~\ref{fig:79}.
For glass bottles tightly packed in an empty box,
packings with rattlers should definitely be avoided
even though the box area was minimal!

The case 
$n = 79$ and several others are marked in Table~\ref{tab:54t213} with two stars
to signify that the packing not only
{\em may} contain a monovacancy, 
but in fact, it {\em must}:
any hexagonal packing of $n$
circles without a monovacancy is provably worse
than the one represented in the table
when $n$ is marked with two stars.

All such cases in the table also happen 
to have $h_{-} = 0$.
Also, a double-star marked entry in Table~\ref{tab:54t213}
always has a discrepancy
between the $n$ computed from the parameters of the packing
structure using \eqref{ncircles},(here it would be
$n=wh$), and actual $n$.
In all the other cases,
the $n$ computed
from the packing structure parameters using \eqref{ncircles}
will always match 
the correct $n$
because there are no
monovacancies.
However, no entry without a monovacancy
is possible in the double-star marked cases.
Another observation:
an entry $n+1$ that follows 
an entry $n$ marked by one or two stars,
always
has the same $w$ and $h$,
and hence its packing fits in the same rectangle 
as that of the entry $n$.

\begin{table}
\begin{center}
\fbox{
\begin{tabular}{r|r|r|r||r|r|r|r||r|r|r|r||r|r|r|r}
$n$&$w$&$h$&$h_{-}$&$n$&$w$&$h$&$h_{-}$&$n$&$w$&$h$&$h_{-}$&$n$&$w$&$h$&$h_{-}$ \\ \hline
54 &14 & 4 & 2     & 94&24 & 4 & 2    &134&34 & 4 & 2    &174& 29 & 6  & 0 \\ 
55 &11 & 5 & 0     & 95&14 & 7 & 3    &135&23 & 6 & 3    &175& 25 & 7  & 0 \\
56 &19 & 3 & 1     & 96&16 & 6 & 0    &136&17 & 8 & 0    &176& 20 & 9  & 4 \\
57 &19 & 3 & 0     &*97&20 & 5 & 3    &137&20 & 7 & 3    &177& 30 & 6  & 3 \\
58 &12 & 5 & 2     & 98&20 & 5 & 2    &138&28 & 5 & 2    &178& 36 & 5  & 2 \\
59 &20 & 3 & 1     & 99&17 & 6 & 3   &*139&16 & 9 & 5    &179& 26 & 7  & 3 \\
60 & 9 & 7 & 3     &100&20 & 5 & 0    &140&16 & 9 & 4    &180& 23 & 8  & 4 \\
*61&21 & 3 & 2     &101&34 & 3 & 1    &141&24 & 6 & 3  &**181& 26 & 7  & 0 \\
62 &21 & 3 & 1     &102&15 & 7 & 3   &*142&29 & 5 & 3    &182& 26 &  7 & 0 \\
63 &13 & 5 & 2     &103&21 & 5 & 2    &143&29 & 5 & 2    &183& 37 &  5 & 2 \\
64 &16 & 4 & 0     &104&12 & 9 & 4    &144&21 & 7 & 3    &184& 23 &  8 & 0 \\
65 &22 & 3 & 1     &105&18 & 6 & 3    &145&29 & 5 & 0    &185& 21 &  9 & 4 \\
66 &17 & 4 & 2     &106&27 & 4 & 2    &146&37 & 4 & 2    &186& 27 &  7 & 3 \\
67 &10 & 7 & 3    &*107&22 & 5 & 3    &147&21 & 7 & 0    &187& 17 & 11 & 0 \\
68 &14 & 5 & 2     &108&22 & 5 & 2    &148&30 & 5 & 2    &188& 24 &  8 & 4 \\
69 &12 & 6 & 3     &109&16 & 7 & 3    &149&17 & 9 & 4    &189& 27 &  7 & 0 \\
70 &14 & 5 & 0     &110&22 & 5 & 3    &150&25 & 6 & 0    &190& 19 & 10 & 0 \\
71 &24 & 3 & 1     &111&19 & 6 & 3    &151&22 & 7 & 3  &**191& 24 &  8 & 0 \\
72 &18 & 4 & 0     &112&16 & 7 & 0    &152&19 & 8 & 0    &192& 24 &  8 & 0 \\
73 &15 & 5 & 2     &113&23 & 5 & 2    &153&31 & 5 & 2    &193& 28 &  7 & 3 \\
74 &11 & 7 & 3     &114&19 & 6 & 0    &154&22 & 7 & 0    &194& 22 &  9 & 3 \\
75 &15 & 5 & 0     &115&23 & 5 & 0    &155&31 & 5 & 0    &195& 20 & 10 & 5 \\
76 &19 & 4 & 0     &116&17 & 7 & 3    &156&20 & 8 & 4    &196& 25 &  8 & 4 \\
77 &26 & 3 & 1     &117&20 & 6 & 3   &*157&23 & 7 & 4  &**197& 22 &  9 & 0 \\
78 &16 & 5 & 2     &118&24 & 5 & 2    &158&23 & 7 & 3    &198& 22 &  9 & 0 \\
**79 &16 & 5 & 0   &119&17 & 7 & 0    &159&27 & 6 & 3   &*199& 29 &  7 & 4 \\
80 &16 & 5 & 0     &120&20 & 6 & 0    &160&20 & 8 & 0    &200& 29 &  7 & 3 \\
81 &12 & 7 & 3    &*121&14 & 9 & 5    &161&23 & 7 & 0    &201& 34 &  6 & 3 \\
82 &21 & 4 & 2     &122&14 & 9 & 4    &162&27 & 6 & 0    &202& 16 & 13 & 6 \\
83 &17 & 5 & 2     &123&18 & 7 & 3    &163&33 & 5 & 2    &203& 23 &  9 & 4 \\
84 &14 & 6 & 0     &124&16 & 8 & 4    &164&21 & 8 & 4    &204& 19 & 11 & 5 \\
85 &17 & 5 & 0     &125&25 & 5 & 0    &165&24 & 7 & 3    &205& 21 & 10 & 5 \\
86 &22 & 4 & 2     &126&21 & 6 & 0   &*166&19 & 9 & 5   &*206& 30 &  7 & 4 \\
87 &15 & 6 & 3     &127&12 & 11& 5    &167&19 & 9 & 4    &207& 30 &  7 & 3 \\
88 &18 & 5 & 2     &128&26 & 5 & 2    &168&34 & 5 & 2    &208& 26 &  8 & 0 \\
89 &30 & 3 & 1     &129&22 & 6 & 3    &169&13 &13 & 0    &209& 19 & 11 & 0 \\
90 &18 & 5 & 0     &130&19 & 7 & 3    &170&17 &10 & 0    &210& 30 &  7 & 0 \\
91 &13 & 7 & 0     &131&15 & 9 & 4    &171&16 &11 & 5   &*211& 24 &  9 & 5 \\
92 &23 & 4 & 0     &132&22 & 6 & 0    &172&25 &7  & 3    &212& 24 &  9 & 4 \\
93 &19 & 5 & 2     &133&27 & 5 & 2    &173&35 &5  & 2    &213& 36 &  6 & 3 

\end{tabular}
}
\caption{
Packings of $n$ circles in a rectangle
for $54 \le n \le 213$ with $w$ circles in a row,
$h$ hexagonal rows,
$h_{-}$ of which
consist of $w-1$ circles each.
Except for the cases marked with stars, 
all packings are the best we could find.
}
\label{tab:54t213}
\end{center}

\end{table}

\begin{figure}
\centering
\includegraphics*[width=7in,height=7in]{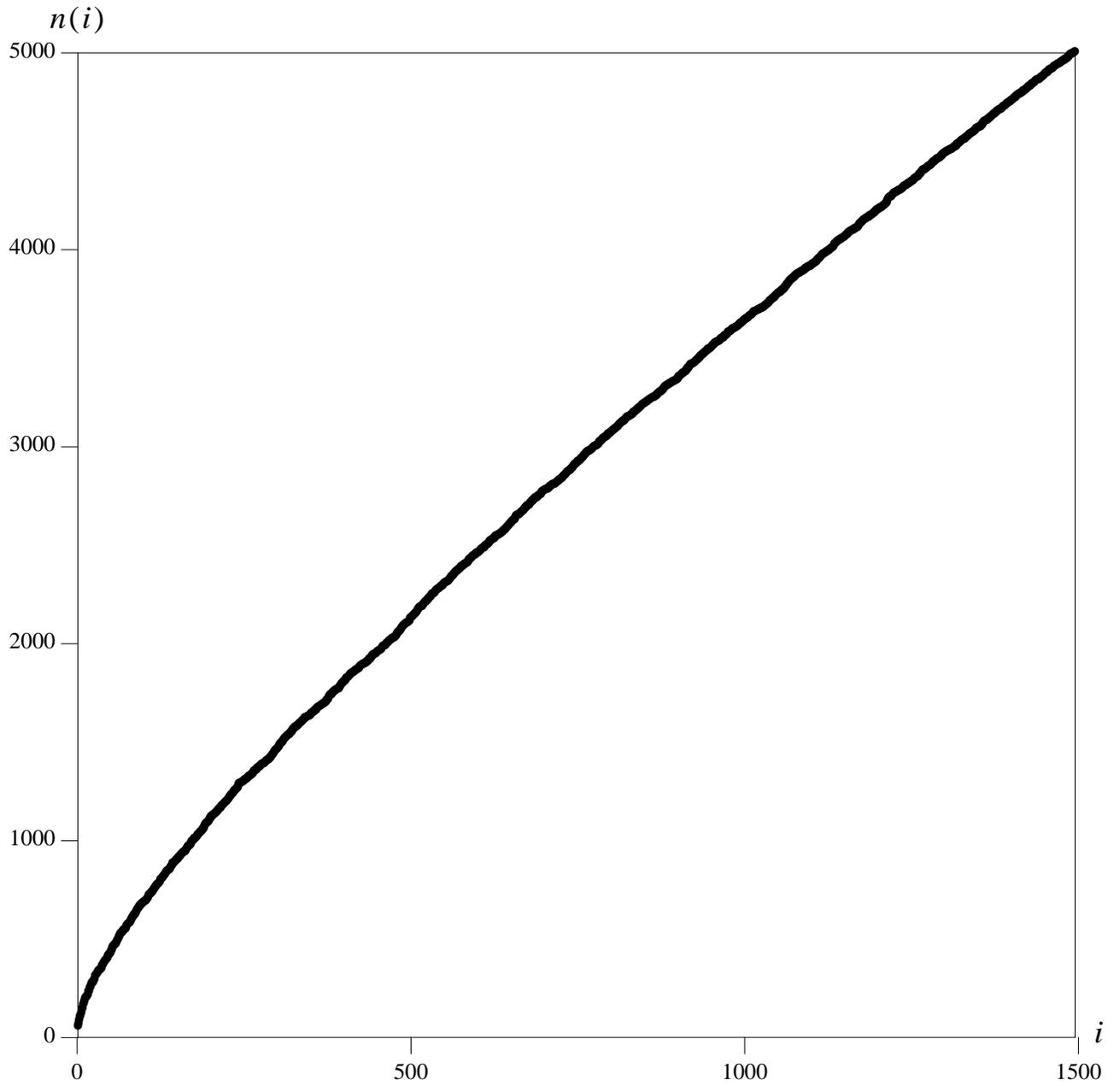}
\caption{Numbers of circles $n(i)$ for which the optimum packing
in a rectangle is not perfectly hexagonal
for $49 \le n(i) \le 5000$ plotted versus the rank $i$.
}
\label{fig:tyt}
\end{figure}

\section{Best packings found for larger $n$}
\hspace*{\parindent} The exhaustive search procedure
of Section~\ref{sec:exhaust}
produces the best packings in the class defined in 
Section~\ref{sec:exhaust}
for values
of $n$ on the order of several thousands.
The features we observed for $n \le 213$ hold
until $n = 317$.
Namely, only one best size rectangle exists for each $n$.
Most of the best packings in the search class
are perfectly hexagonal and unique for their $n$.
Each such packing can be
described by the set $w$, $h$, and $h_{-}$,
with the latter parameter taking on up to three possible
values:
\begin{equation}
\label{hminus}
\begin{array}{lllll}
   \mbox{if $h$ is even} &\mbox{ then }
       &h_{-} =  0 ~~\mbox{ or }&h_{-} = h/2 &\\
   \mbox{if $h$ is odd}  &\mbox{ then } &h_{-} =  0 ~~
        \mbox{ or }&h_{-} = \lfloor h/2 \rfloor ~~
        \mbox{ or }&h_{-} = \lfloor h/2 \rfloor + 1
\end{array}
\end{equation}
A few exceptional cases have a single monovacancy.
These are
similar to the one-star marked cases with odd $h$
and with two options for $h_{-}$ 
(the choice of the option defines the presence or 
absence of the vacancy),
or they are similar to the two-star marked cases with $h_{-} = 0$.

The case $n = 317$ deviates from this pattern.
Here, the best packing in the class has parameters
$w = 27$, $h = 12$, $h_{-} = 6$, and $d = 1$.
Unlike the one-star marked cases, $h$ is even,
and unlike the two-star marked cases, $h_{-}$ is non-zero.
Such a new type of packing recurs with one hole for
$n = 334$ ($w = 34$, $h = 10$, and $h_{-} = 5$)
and then 
for 
$n = 393$ ($w = 40$, $h = 10$, and $h_{-} = 5$),
but in the latter case with {\em two} monovacancies.
In other words, the best possible packing
for 393 circles in a rectangle in the class of hexagonal
packings with possible monovacancies must have
two of them.

With two holes,
we have more freedom to improve the quality of the
packing than having just one hole.
Several ways are possible, including combining the ones
described before for cases of smaller $n$. 
Same as for the
other packings with monovacancies,
the optimal packing for $n = 393$ is not known,
but the exhaustive search
in our experiments proves that it cannot be
a hexagonal packing.

Naturally, the best packing (in the search class) of $n = 394$ circles
also has parameters $w = 40$, $h=10$, $h_{-} = 5$ and $d = 1$,
i.e., a single monovacancy.
The same rectangle also accommodates the best packing in the class 
of $n = 395$ circles, with the same $w = 40$, $h=10$, $h_{-} = 5$.
Since the latter packing is without holes, 
it is conceivable that it is optimal.

The next outstanding case is $n = 411$
where $w = 38$, $h = 11$, $h_{-} = 6$, and $d = 1$.
This is the smallest $n$ 
where the $h$ is odd and $h_{-} = \lfloor h/2 \rfloor + 1$,
that is, the packing looks like the one in Figure~\ref{fig:49}a,
but unlike the latter it has a hole.
An equivalent packing that looks like the one in Figure~\ref{fig:49}b,
where $h_{-} = \lfloor h/2 \rfloor = 5$ also exists, 
but unlike the latter it has not one but 
two holes.
For next value $n = 412$, we have a standard one-star marked situation
with the same sizes of the rectangle, i.e.,
$w = 38$, $h = 11$ and 
$h_{-} = \lfloor h/2 \rfloor = 5$
for the variant with one hole 
and $h_{-} = \lfloor h/2 \rfloor + 1 = 6$
for the variant without holes.

The smallest $n$ for which as many as three monovacancies exist
in the best packing in the class is $n = 717$.
The corresponding packing parameters are
$w = 48$, $h = 15$, $h_{-} = 0$ and no hexagonal packing
with a smaller number of holes can be better or even as good as this
packing.
Similarly,
with four monovacancies,
the smallest $n$ 
is
$n = 2732$ 
($w =  86$, $h = 32$, and $h_{-} = 16$),
and for five monovacancies, $n = 2776$ is the smallest
($w = 103$, $h = 27$, and $h_{-} = 0$).
No optimum packing in the class for $n \le 5000$ has six or more
monovacancies.

In all, we found 1495 values of $n$ on the interval
$1 \le n \le 5000$ for which the best packing in the class must
or can have monovacancies. 
As explained above, for each such $n$ the best packing  provably
cannot be of a pure square-grid or hexagonal pattern.
It is not proven though, but we believe
it to be also true, that those 1495 values of $n$ are {\em all} such
irregular values among the considered 5000 values.

The chance to encounter such an irregular $n$
seems to increase with $n$.
For example, here are all the experimentally found 
irregular values of $n$ among 100 consecutive values on the intervals
$401+1000k \le n \le 500 +1000k$ for $k=0,1,2,3,4$:

17 values for $401 \le n \le 500$:
 
\noindent 
409, 411, 412, 421, 422, 
433, 439, 453, 454, 461, 
463, 467, 471, 478, 487, 
489, 499.
  
24 values for $1401 \le n \le 1500$:

\noindent 
1401, 1402, 1405, 1409, 1412, 
1414, 1423, 1427, 1429, 1434,
1446, 1447, 1451, 1453, 1457,
1459, 1466, 1468, 1477, 1483,
1486, 1487, 1489, 1497. 
  
33 values for $2401 \le n \le 2500$:

\noindent 
2401, 2402, 2406, 2411, 2419,
2421, 2423, 2428, 2429, 2435,
2437, 2439, 2441, 2443, 2446,
2452, 2454, 2455, 2456, 2458,
2462, 2467, 2469, 2474, 2476,
2477, 2479, 2481, 2487, 2491,
2493, 2495, 2497.
  
33 values for $3401 \le n \le 3500$:

\noindent 
3407, 3409, 3411, 3412, 3414,
3415, 3418, 3421, 3425, 3428,
3431, 3433, 3436, 3442, 3446
3447, 3453, 3455, 3459, 3461,
3464, 3467, 3469, 3473, 3476, 
3479, 3481, 3487, 3489, 3490,
3493, 3494, 3499.
  
38 values for $4401 \le n \le 4500$: 

\noindent 
4401, 4404, 4405, 4409, 4411,
4414, 4417, 4419, 4421, 4426,
4430, 4434, 4436, 4438, 4441,
4443, 4447, 4450, 4453, 4456,
4457, 4458, 4461, 4462, 4467,
4468, 4474, 4476, 4479, 4483,
4486, 4487, 4491, 4492, 4493,
4495, 4497, 4499.

Figure~\ref{fig:tyt} represents these 1495 irregular values $n(i)$,
$1 \le n(i) \le 5000$,
beginning with $n(1)=49$ and ending with $n(1495)=4999$, 
by plotting points 
\{abscissa$~=~i$,\ ordinate$~=~n(i)$\}.
\section{The optimum aspect ratio for large $n$}
\begin{figure}
\centering
\includegraphics*[width=6.3in,height=8.2in]{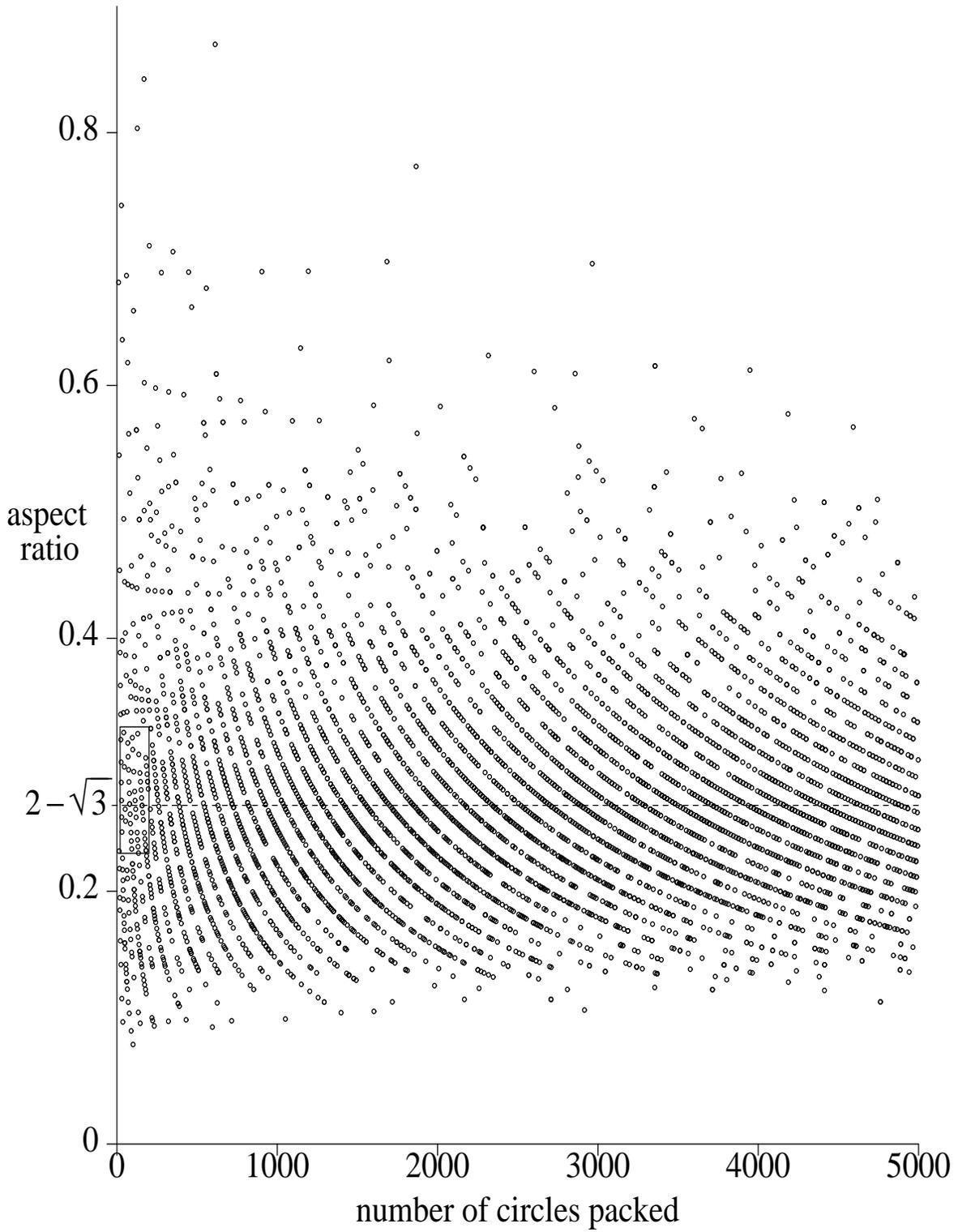}
\caption{Values of the aspect ratio for optimal hexagonal packings.}
\label{fig:hovbst}
\end{figure}

\hspace*{\parindent} 
Figure~\ref{fig:hovbst}
contains for each $n$,
a data point with coordinates
\\
\\
$~~~~~~~~~~~~~~$(abscissa = $n$, ordinate = best aspect ratio found for this $n$).
\\
\\
All $n \le 5000$ are represented with the exception
of a few $n$ where the best packings
found are of the square grid type.
Also, the hybrid packings like those in Figure~\ref{fig:15}b
are excluded. 
In this way we assure the uniqueness of the
best aspect ratio for each $n$ represented.
The aspect ratios for the best packings in the
search class are described in Section~\ref{sec:exhaust}.
The changes to the aspect ratios are due to the
$\delta$ shrinkages of the rectangles
in those cases with monovacancies, such as
those in Figures~\ref{fig:49}, \ref{fig:79}, and \ref{fig:5h3}.
Presumably these are small in number and should not be noticeable
in Figure~\ref{fig:hovbst}.

\begin{figure}
\centering
\includegraphics*[width=7.3in,height=6.2in]{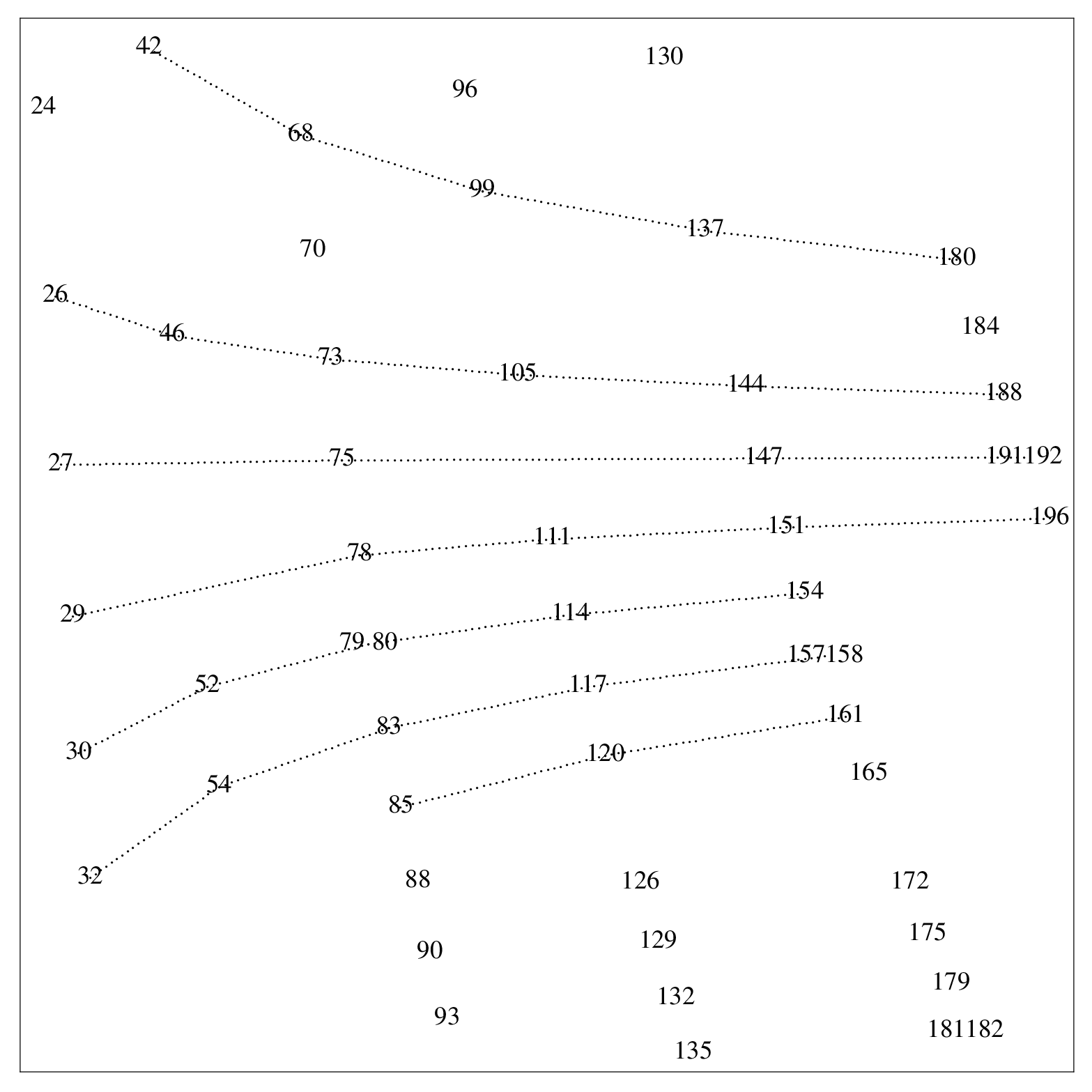}
\caption{A box selected in Figure~\ref{fig:hovbst} magnified.
Each data point in Figure~\ref{fig:hovbst} is replaced
here with the corresponding number of circles.}
\label{fig:boxa}
\end{figure}

The data points tend to form patterns 
of
descending and ascending ``threads.''
To examine the threads in detail,
a rectangular box which is close to the $y$-axis 
in Figure~\ref{fig:hovbst} is magnified
and represented in
Figure~\ref{fig:boxa}.
All data points present in the box can be found in 
Tables~\ref{tab:1t53} and \ref{tab:54t213}.

The steep downward threads
correspond to configurations with fixed 
height, and with widths increasing by 1 from 
one point to the next as we move down the thread.
Each such thread eventually terminates in either direction
which means that
the optimum rectangle cannot be too 
flat or too tall.
The less steep upward and downward threads (they are dotted)
correspond to sequences ($w,h$), where from one point to the next
$w$ increases by 3 or 6 and $h$ increases by 1 or 2, respectively.

\begin{figure}
\centering
\includegraphics*[width=5in,height=3.9in]{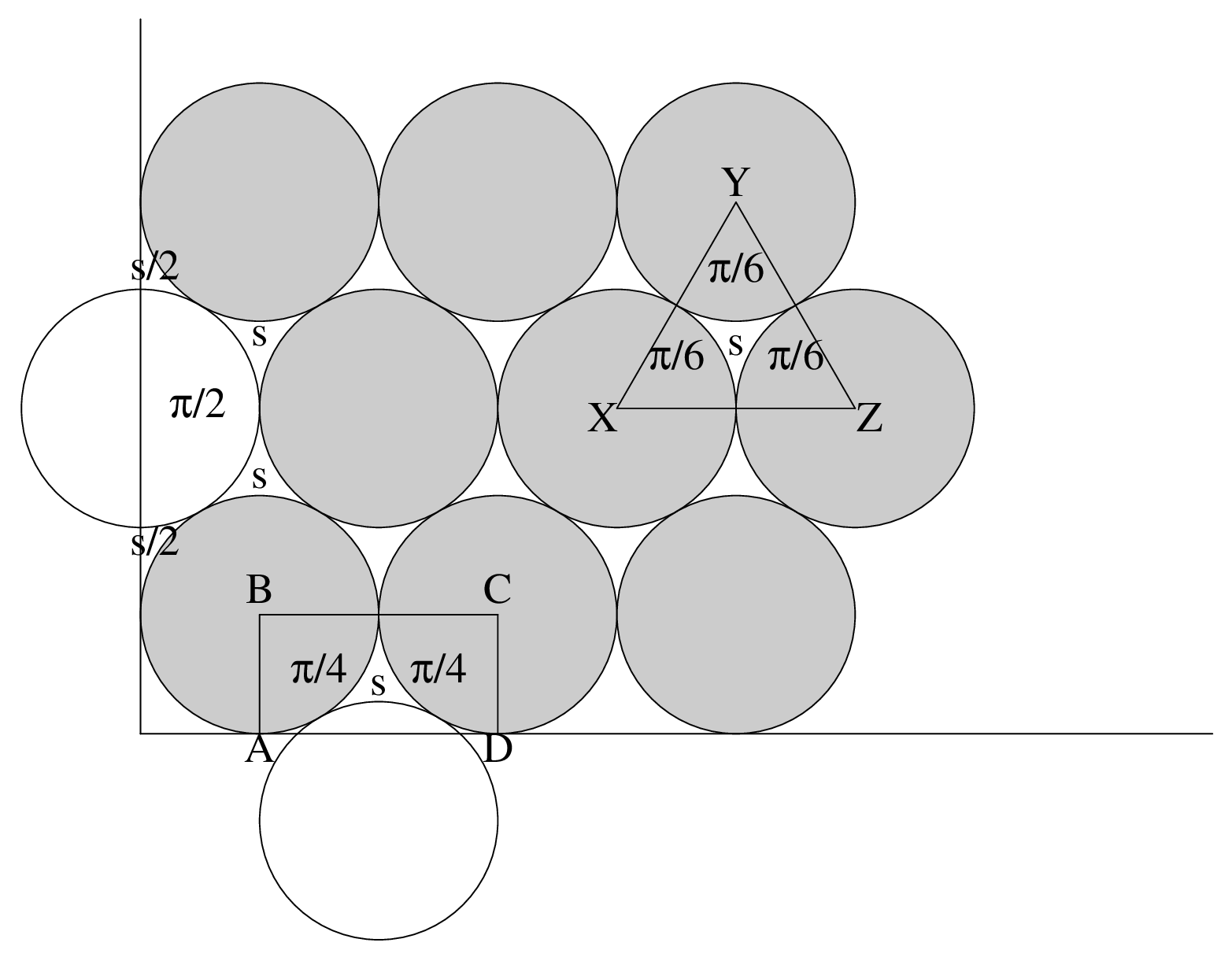}
\caption{Calculation of the waste area adjacent 
to the sides of the rectangle.}
\label{fig:waste}
\end{figure}

To compute the optimum aspect
ratio of the rectangle that encloses the optimum packing for 
large $n$, we
analyze the area wasted along the rectangle sides.
Note that in the infinite hexagonal packing,
the uncovered area 
is $s = \sqrt 3 - \pi /2$ per each $\pi /2$ of the
covered area.
This is obvious from examining a single triangle $XYZ$ 
in Figure~\ref{fig:waste} and observing that
the entire infinite packing is composed
of such triangles.
The waste $s$ here is
the area of the central triangle formed by
three circular arcs and it is
equal to
the full area $\sqrt{3}$ of triangle $XYZ$
(assuming circle radius is 1) minus $\pi /2$, which is
the total area of the three $\pi /6$ sectors covered.
(By the way, we remind the reader that the density of the infinite
hexagonal packing is 
$(\sqrt 3 - s) / (\pi/2) = \pi / (2 \sqrt 3 )  = 0.90689968...$.)

The structure of the uncovered area in the infinite hexagonal packing
can be also understood if we think of
each covered circle 
bringing with it two curved uncovered triangles $s$,
those adjacent to this circle on its right-hand side.
With such bookkeeping, 
each triangle $s$ will be counted exactly once.
We conjecture that these two triangles $s$ per each circle
is the unavoidable, 
i.e., fixed, waste in any finite hexagonal packing
carved out by a rectangle.
There will be additional, i.e., variable,
waste along the rectangle sides
and we are trying to minimize this variable waste.

Along the bottom and top sides, for each 2 units of length,
like the side $AD$ of
rectangle $ABCD$ in Figure~\ref{fig:waste},
the additional waste is the area of rectangle
$ABCD$ minus two covered quarter-circle areas and minus
$s$. 
This $s$ represents one of the two curved triangles attached to the
right-hand side of the circle with the center at $B$ 
and as an unavoidable waste should not be included.
Hence the waste per unit length of the top and bottom sides
is $a = (2-\pi /2 -s)/2 = (2 - \sqrt 3 )/2$.

The waste along the left- or right-hand sides per two alternating rows,
i.e., per $2 \sqrt{3}$ units of length,
is the area of the semicircle that is cut in half by 
the side 
plus or minus additional area depending on the side.
Along the left-hand side,
all additional uncovered area is 
additional waste,
since we assume that the curved triangles $s$ are attached
to the right-hand side of the covered circles.
The addition consists of two triangles $s$ and two halves
of triangles $s/2$ as shown in Figure~\ref{fig:waste}.
The left-hand side waste per $2 \sqrt{3}$ of length is thus
$\pi/2 + 3s$.
Along the right-hand side we subtract from the area of
the semicircle the outstanding  two half triangles $s/2$
because these are necessary in the infinite packing.
The right-hand side waste per $2 \sqrt{3}$ units of length is thus 
$\pi/2 -s$.
Averaging this for both sides and dividing by $2 \sqrt{3}$,
we have $b = ( \pi/2 + s)/ (2 \sqrt 3 ) = 1/2$ as the additional
waste per unit length of left-hand or right-hand sides.

Now we should take $a/b = 2 - \sqrt{3}$ as the ratio
of height over width that yields the optimal balance 
between the waste along the sides of the enclosing
rectangle and hence the minimum area.
Figure~\ref{fig:hovbst} includes a horizontal
line at the level of the aspect ratio
$2 - \sqrt{3}$. 

\section{Concluding remarks}
As pointed out in the beginning of the paper, we hope that
these computational results will lead to proofs of 
optimality for larger (or even infinite classes of) $n$.

Our experiments 
suggest
that the frequency of occurrence of non-hexagonal 
best packings increases with $n$.
It is not clear whether or not the frequency
has a limit as $n \rightarrow \infty$ and if it has,
whether of not this limit is smaller than 1.
While it is not easy to find
small $n \ge 14$ for which the best packing is not perfectly hexagonal,
it might be more difficult to find large $n$
for which the best packing {\em is} perfectly hexagonal.
Do there exist infinitely many $n$,
for which the densest packing of $n$ circles in a rectangle
is hexagonal?

A related phenomenon is 
conjectured to hold for $n = N(k)$ of the form $N(k)=\frac{1}{2}(a_k+1)(b_k+1)$
where $a_k$ and $b_k$ are given by:

\[ a_1 = 1, a_2 = 3, a_{k+2} = 4a_{k+1} - a_k\]
\[ b_1 = 1, b_2 = 5, b_{k+2} = 4b_{k+1} - b_k\]

\noindent 
so that, $N(2) = 12, N(3) = 120, N(4) = 1512$, etc. The fractions
$\frac{a_k}{b_k}$ are actually (alternate) convergents to
$\frac{1}{\sqrt 3}$, and it has been conjectured by Nurmela et al.
\cite{NOR} that for these $n$, a ``nearly'' hexagonal packing of
$n$ circles in a square they describe is in fact optimal. 

In our case, one should seek alternate  convergents to
$\sqrt 3  + 3/2$
which yields sequences $a_k$ and $b_k$ given by:

\[a_k = 2v_{k+1} - v_k,\ \   b_k = 2v_k,\]
\noindent
where
\[v_0 = 0,\ \ v_1 = 1,\ \ v_{k+2} = 4v_{k+1} - v_k,\ \ k = 2, 3,... \]
Thus, 
 \[ (a_1,b_1) = (7,2),\ \ (a_2, b_2) = (26,8),\ \  (a_3, b_3) = (97,30) ....\]
\noindent
so that $N(2) = 208,\ \ N(3) = 2910$, etc. 
No such $N(k)$ for $k = 2,3,...$ is indeed found to be irregular
in our experiments.
The best packing found experimentally for such an $N(k)$
has $h = b_k$ alternating rows of full length 
$w = a_k$ with $h_{-} = 0$, $s = 0$,
in the notation of 
Section~\ref{sec:exhaust} and 
Tables~\ref{tab:1t53} and ~\ref{tab:54t213}.


\end{document}